\documentclass{amsart}
\usepackage{amscd,amssymb,amsopn,amsmath,amsthm,graphics,amsfonts,accents,enumerate,verbatim,calc}
\usepackage[dvips]{graphicx}
\usepackage{pgfplots}
\usepackage[colorlinks=true,linkcolor=red,citecolor=blue]{hyperref}
\usepackage[all]{xy}
\usepackage{cite}
\usepackage{tikz}
\usepackage{tikz-cd}
\usepackage{mathrsfs}

\calclayout

\newcommand{\rt}{\rightarrow}
\newcommand{\lrt}{\longrightarrow}

\newcommand{\st}{\stackrel}

\newcommand{\La}{\Lambda}

\newcommand{\SA}{\mathscr{A}}
\newcommand{\SB}{\mathscr{B}}
\newcommand{\SC}{\mathscr{C}}
\newcommand{\SD}{\mathscr{D}}
\newcommand{\SE}{\mathscr{E}}
\newcommand{\SF}{\mathscr{F}}

\newcommand{\SR}{\mathscr{R}}

\newcommand{\ST}{\mathscr{T}}

\newcommand{\SX}{\mathscr{X}}

\newcommand{\CE}{\mathcal{E}}

\newcommand{\CR}{\mathcal{R} }

\newcommand{\CT}{\mathcal{T} }

\newcommand{\Mod}{{\rm{Mod\mbox{-}}}}

\newcommand{\mmod}{{\rm{{mod\mbox{-}}}}}

\newcommand{\Inj}{{\rm{Inj}\mbox{-}}}
\newcommand{\Prj}{{\rm{Prj}\mbox{-}}}
\newcommand{\prj}{{\rm{prj}\mbox{-}}}

\newcommand{\Prod}{{\rm{Prod}}}

\newcommand{\add}{{\rm{add}\mbox{-}}}

\newcommand{\Add}{{\rm{Add}\mbox{-}}}

\newcommand{\Fac}{{\rm{Fac}}}

\newcommand{\Sperp}{S^{\mbox{\small{\rm{perp}}}}}

\newcommand{\Pres}{{\rm Pres}}

\newcommand{\pd}{{\rm{pd}}}

\newcommand{\Coker}{{\rm{Coker}}}
\newcommand{\Ker}{{\rm{Ker}}}

\newcommand{\Hom}{{\rm{Hom}}}
\newcommand{\Ext}{{\rm{Ext}}}

\newcommand{\Gen}{{\rm{Gen}}}
\newcommand{\Cogen}{{\rm{Cogen}}}

\newcommand{\Copres}{{\rm Copres}}

\theoremstyle{plain}
\newtheorem{theorem}{Theorem}[section]
\newtheorem{corollary}[theorem]{Corollary}
\newtheorem{lemma}[theorem]{Lemma}

\newtheorem{proposition}[theorem]{Proposition}

\theoremstyle{definition}
\newtheorem{definition}[theorem]{Definition}

\newtheorem{remark}[theorem]{Remark}
\newtheorem{setup}[theorem]{Setup}

\theoremstyle{plain}

\theoremstyle{definition}

\numberwithin{equation}{section}

\begin{document}
\title[Tilting subcategories and silting modules]{Extending ($\tau$-)tilting subcategories \\ and (co)silting modules}

\author[J. Asadollahi, F. Padashnik, S. Sadeghi and H. Treffinger]{J. Asadollahi, F. Padashnik, S. Sadeghi and H. Treffinger}

\address{Department of Pure Mathematics, Faculty of Mathematics and Statistics, University of Isfahan, P.O.Box: 81746-73441, Isfahan, Iran}
\email{asadollahi@sci.ui.ac.ir, asadollahi@ipm.ir}
\email{f.padashnik@sci.ui.ac.ir}
\email{so.sadeghi@sci.ui.ac.ir }

\address{Universit\'{e} de Paris, B\^{a}timent Sophie Germain 5, rue Thomas Mann 75205, Paris Cedex 13, FRANCE}
\email{treffinger@imj-prg.fr}

\makeatletter \@namedef{subjclassname@2020}{\textup{2020} Mathematics Subject Classification} \makeatother



\begin{abstract}
Let $B$ be a finite dimensional algebra and $A=B[P_0]$ be the one-point extension algebra of $B$ with respect to the finitely generated projective $B$-module $P_0$. The categories of $B$-modules and $A$-modules are related by two adjoint functors $\CR$ and $\CE$, called the restriction and the extension functors, respectively. Based on the nice homological properties of these two functors, restriction and extension of some notions such as tilting and $\tau$-tilting modules have been studied in the category of finitely presented modules, i.e. small mod. In this paper, we investigate the behaviour of tilting and support $\tau$-tilting subcategories with respect to these two functors. Moreover, we investigate the restriction and the extension of special related modules such as finendo quasi-tilting modules, silting modules, and cosilting modules. Our studies will be done in the category of all modules, which will be called large Mod. Based on such study, in addition to the new results, classical results are extended, not only from modules to subcategories but also from small mod to large Mod.
\end{abstract}

\maketitle


\section{Introduction}
All algebras in this paper are assumed to be finite dimensional over an algebraically closed field $k$. Let $B$ be such an algebra and $A=B[P_0]$ be the one-point extension algebra of $B$ with respect to a fixed finitely generated projective $B$-module $P_0$. It is shown by Assem, Happel and Trepode \cite{AHT} that one can construct a tilting $B$-module, resp. a tilting $A$-module, starting from a basic tilting $A$-module, resp. a basic tilting $B$-module. To do this they considered two functors $\SR: \mmod A \lrt \mmod B$ and $\SE: \mmod B \lrt \mmod A$, called restriction and extension functors, respectively. It turned out that these two functors, being an adjoint pair, have nice homological properties. P. Suarez \cite{Su} followed this approach and studied the restriction and extension of (support) $\tau$-tilting modules. The notion of $\tau$-tilting theory was introduced by Adachi, Iyama and Reiten in \cite{AIR} as a new approach for  studying two classical branches of the representation theory of finite dimensional algebras, namely tilting theory and Auslander-Reiten theory.

Based on the importance of the notions of tilting and $\tau$-tilting modules, they have been studied, and also generalized, in different settings.

For instance, both notions have been extended to the categorical level. Beligiannis \cite{B} developed a fully general tilting theory in an arbitrary abelian category in an extensive manuscript that is started in 2004 and announced at several conferences later but still not made publicly available, see \cite[Remark 1.9]{BBOS}. More recently, in \cite{BBOS}, the authors reconsidered tilting theory in an abelian category with enough projective objects. On the other hand, the notion of $\tau$-tilting subcategories were studied in \cite{IJY} for functor categories. The more general notion of $\tau$-tilting subcategories in an abelian category is studied in \cite{LZh}. See preliminaries section for more details.

Furthermore, silting modules, introduced and studied in \cite{AMV}, provide a wide generalisation of tilting modules as well as support $\tau$-tilting modules over finite dimensional algebras to arbitrary rings. By Proposition 3.10 of \cite{AMV} all tilting modules are silting,  and by Proposition 3.15 of \cite{AMV} over an artin algebra, a finite dimensional module is silting if and only if it is support $\tau$-tilting. Moreover, in the same paper, an extension of silting modules, i.e. finendo quasi-tilting modules, is introduced and studied. While by \cite[Proposition 3.10]{AMV} all silting modules are finendo quasi-tilting, it is known that the inclusion of silting modules in the class of finendo quasi-tilting modules is proper \cite[Example 5.4]{AH}. As the categorical dual of silting modules, the authors of \cite{BP} introduced the notion of cosilting modules.
We remark that there is a notion of cofinendo quasi-cotilting modules in the literature; by \cite{ZW} all quasi-cotilting modules are cofinendo and by Theorem 4.18 of \cite{ZW1}, quasi-cotilting modules and cosilting modules are the same.

Note that Keller and Vossieck in \cite{KV} introduced silting objects in triangulated categories, as important tools in the study of homotopy or derived categories
and showed that they are in correspondence with other concepts such as (co-)t-structures or simple-minded collections of objects.

A natural attempt now is to investigate the behaviour of tilting and $\tau$-tilting subcategories on one hand and silting modules, finendo quasi-tilting modules, and cosilting modules on the other hand, under the restriction and extension functors. This will be the main theme of this paper. To be able to investigate such behaviour, we need to work within the category of all modules, say large Mod.

The paper is structured as follows. Section \ref{Sec 2-One point} is devoted to study the module category of one-point extension algebras. We show that similar to the nice homological properties of the restriction and extension functors hold true in the large module category. In particular, when $A=B[P_0]$ is the one point extension algebra of $B$ with respect to the finitely generated projective $B$-module $P_0$, then we have the following recollement
\begin{equation}
\xymatrix{\Mod k \ar[rrr]^{i_*}  &&& \Mod A \ar[rrr]^{\textcolor{black}{\mathcal{R}}} \ar@/^2pc/[lll]^{v = \Hom_A(k,-)} \ar@/_2pc/[lll]_{u = k\otimes_A-} &&& \Mod B. \ar@/^2pc/[lll]^{\textcolor{black}{\mathcal{E}}} \ar@/_2pc/[lll]_{\mathcal{L}=Ae_B{\otimes}_B -}}
\end{equation}
In Section \ref{Sec 3-Tilting Sub} we study tilting subcategories. It is shown that in an abelian category with enough projective objects, the definition of tilting subcategories introduced in \cite{BBOS} is equivalent to the definition of \cite{B}. Using this fact, in Section \ref{Sec 4-TilSub One point} the behaviour of tilting subcategories with respect to the restriction and extension functor will be studied. Results of this section, provide an extensive generalization of the main result of \cite{AHT}.

Then we turn our attention to the study of restriction and extension of certain important modules in the module category of an algebra. The modules we consider are silting modules, finendo quasi-tilting modules and cosilting modules. This will be done in Section \ref{Sec 5-Special Mod}.

$\tau$-tilting subcategories of abelian categories are a wide generalization of $\tau$-tilting modules. Let $\SA$ be an abelian category with enough projective objects. By Definition 2.1 of \cite{LZh} an additively closed full subcategory $\ST$ of $\SA$ is called a weak support $\tau$-tilting subcategory if $\Ext^1_{\SA}(T_1,\Fac(T_2))=0,$ for all $T_1, T_2 \in \ST$ and for every projective object $P$ in $\SA$, there exists an exact sequence $P \st{f}{\lrt} T^0 \lrt T^1 \lrt 0$ such that $T^0$ and $T^1$ are in $\ST$ and $f$ is a left $\ST$-approximation of $P$. In Section \ref{Sec 6-tau-tilSub} the restriction and extension of $\tau$-tilting subcategories will be studied.
On the other hand, there is a bijection between the collection of all equivalence classes of $\tau$-tilting subcategories of $\Mod R$ and the collection of all equivalence classes of finendo quasi-tilting $R$-modules, when $R$ is an arbitrary ring \cite[Theorem 8.1.10]{AST}. In view of this result, we are able to show that restriction of a finendo quasi-tilting module $T$ in $\Mod A$, i.e. $\CR T$, is a finendo quasi-tilting module in $\Mod B$.

By Theorem 2.29 of \cite{BBOS}, we know that there are bijections between tilting subcategories of $\SA$, where $\SA$ is an abelian category with enough projective objects, and special triples of full subcategories of $\SA$, the so-called cotorsion torsion triples. Moreover,  Theorem 5.7 of \cite{AST} implies that there is a bijection between the collection of all support $\tau$-tilting subcategories of $\SA$ and the collection of all $\tau$-cotorsion torsion triples of $\SA$. Based on these bijections we study explicitly the behaviour of these triples under the restriction and extension functors. This is the content of the last section of the paper.

{\sc Notations and conventions.} Throughout the paper, all rings are associative with identity. Let $R$ be such a ring. The category of all (left) $R$-modules will be denoted by $\Mod R$. We let $\Prj R$, resp. $\Inj R$, denote the full subcategory of $\Mod R$ consisting of all projective, resp. injective, $R$-modules.

For a class $\SX$ in $\Mod R$, let $\Add \SX$, resp. $\add \SX$, denote the class of all modules isomorphic to a direct summand of an arbitrary direct sum, resp. finite direct sum, of copies of modules in $\SX$. Also let $\Gen(\SX)$, resp. $\Fac(\SX)$, be the subcategory of $\Mod R$ consisting of all $R$-modules isomorphic to an epimorphic images of modules in $\Add \SX$, resp. in $\add \SX$. Obviously we have $\Gen(\SX)= \Fac(\Add \SX)$.
Dually, let $\Prod(\SX)$ denote the subcategory of $\Mod R$ consisting of all modules isomorphic to a direct summand of an arbitrary direct product of copies of modules in $\SX$. We also let $\Cogen(\SX)$ to be the subcategory of $\Mod R$ consisting of all $R$-modules that can be embedded into a direct product of copies of modules in $\SX$.

Finally, for an $R$-module $M$, let $\Pres(M)$ denote the subcategory of $\Mod R$ consisting of all modules that admit an $\Add M$-presentation, i.e. all modules $X$, for them there exists an exact sequence $M_1 \lrt M_0 \lrt X \lrt 0$, such that $M_0, M_1 \in \Add M$ and let $\Copres(M)$ be the subcategory of $\Mod R$ consisting of all modules that admit an $\Prod(M)$-copresentation.

\section{One-point extension algebras}\label{Sec 2-One point}
Let $k$ be an algebraically closed field and $B$ be a finite dimensional $k$-algebra. Let $\Mod B$ denote the category of (left) $B$-modules and $\mmod B$ denote its subcategory consisting of finitely presented $B$-modules.

Let $\{e_1, \ldots, e_n\}$ be a basic set of primitive idempotents of $B$. It is known \cite[Theorem 27.11]{AF} that projective $B$-modules are direct sums of (possibly infinitely many) copies of $Be_i$, for $1 \le i\le n$. When $B$ is a bounded path algebra of a quiver $Q$, these idempotents are in bijection with the vertices of $Q$. So, the only difference between the big projectives, i.e. projectives in $\Mod B$, and small projectives, i.e. projectives in $\mmod B$, are the multiplicities, but there are no new indecomposable projectives showing up.

Let $P_0$ be a fixed finitely generated projective $B$-module. The one-point extension of $B$ by $P_0$, which is denoted by $A=B[P_0]$, is the matrix algebra
\[ A=\bordermatrix{& & \cr & B & P_0\cr & 0 & k\cr}\]
with the ordinary addition of matrices and the multiplication induced from the module structure of $P_0$.

It is known that $B$ is a full convex subcategory of $A$ and there exists a unique indecomposable projective $A$-module $\tilde{P}$ which is not a projective $B$-module. Also, the simple top $S$ of $\tilde{P}$ is an injective $A$-module of projective dimension at most one.

Let $e_B$ denote the identity of $B$. The following two functors
\[{\mathcal{R}}=\Hom_A(A e_B, -): \Mod A \lrt \Mod B\]
and
\[ \ \mathcal{E}=\Hom_B(e_B A,-): \Mod B \lrt \Mod A\]
are called restriction and extension functors, respectively.

It is known that $(\mathcal{R}, \CE)$ is an adjoint pair. More generally, they fit into the following recollement
\begin{equation}\label{recollement}
\xymatrix{\Mod k \ar[rrr]^{i_*}  &&& \Mod A \ar[rrr]^{\textcolor{black}{\mathcal{R}}} \ar@/^2pc/[lll]^{v = \Hom_A(k,-)} \ar@/_2pc/[lll]_{u = k\otimes_A-} &&& \Mod B. \ar@/^2pc/[lll]^{\textcolor{black}{\mathcal{E}}} \ar@/_2pc/[lll]_{\mathcal{L}=Ae_B{\otimes}_B -}}
\end{equation}
The functor $\CR$ is exact since it is both a right and a left adjoint. Furthermore, since $e_B A$ is a projective $B$-module, the functor $\CE$ is also exact. It follows form the general properties of the functors appearing in a recollement that $\CR$ preserves projective and injective modules and $\CE$ preserves injective modules. Moreover, we can embed $\Mod B$ in $\Mod A$ under the usual embedding functor. In particular, $\CR X$ is a submodule of $X$.

We follow the convention of \cite{AHT} and use letters $X$, $Y$ and $Z$ to denote the $A$-modules and use letters $M$, $N$ and $L$ to denote the $B$-modules.

\begin{remark}
By restricting the categories of the recollement \ref{recollement} to the corresponding subcategories consisting of finitely presented modules we obtain the following recollement
\begin{equation}\label{f.g.recollement}
\xymatrix{\mmod k \ar[rrr]^{i_*}  &&& \mmod A \ar[rrr]^{\textcolor{black}{\SR=\Hom_A(A e_B, -)}} \ar@/^2pc/[lll]^{ \Hom_A(k,-)} \ar@/_2pc/[lll]_{ k\otimes_A-} &&& \mmod B \ar@/^2pc/[lll]^{\textcolor{black}{\SE=\Hom_B(e_B A,-)}} \ar@/_2pc/[lll]_{Ae_B{\otimes}_B -}}
\end{equation}
which was considered in \cite{AHT} and \cite{Su}.
\end{remark}

The proofs of the following properties of the restriction and extension functors are similar to the proofs of Lemma 2.2 and Proposition 3.2 of \cite{AHT}, so we skip their proofs.

\begin{lemma}\label{2.2}
Let $X$ be an $A$-module and $M$ be a $B$-module. There are short exact sequences
\[0\lrt \CR X\lrt X \lrt S''\lrt 0\]
and
\[ 0\lrt M \lrt \CE M \lrt S'\lrt 0\]
in $\Mod A$, where $S'$ and $S''$ are in $\Add S$. They are called the restriction sequence of $X$ and the extension sequence of $M$, respectively.
\end{lemma}

Following proposition collects some of the basic homological properties of the restriction and extension functors. The proof of $(1)$ can be found in Proposition 3.2, Theorem 3.10,  of \cite{P}. The proof of others are similar to the proof of Corollary 3.5, Proposition 3.6, Corollary 3.7 and Lemma 4.5 of \cite{AHT}. So we skip the proofs.

Recall that the right perpendicular category of $S$ is the full subcategory of $\Mod A$ defined by
\[\Sperp:=\lbrace X\in \Mod A~~\vert~~\Hom_A(S, X)=0, ~\Ext^1_A(S, X)=0\rbrace.\]

\begin{proposition}\label{AHT}
Let $X$ and $Y$ be $A$-modules and $M$ be a $B$-module.
\begin{itemize}
  \item[$(1)$] If $X\in \Sperp$, then there is a functorial isomorphism $X\cong \CE\CR X$. Moreover, for all $j\geq 0$, there is an isomorphism \[\Ext^j_A(\CE M, X)\cong \Ext^j_B(M, \CR X).\]
  \item[$(2)$] For all $j\geq 0$, there is an isomorphism \[ \Ext^j_A(X, \CE M)\cong \Ext^j_B(\CR X, M).\]
  \item[$(3)$] There is an epimorphism \[ \Ext^1_A(X, Y)\lrt \Ext^1_B(\CR X, \CR Y).\]
  \item[$(4)$] For each $j\geq 2$, there is an isomorphism \[\Ext^j_A(X, Y)\cong \Ext^j_B(\CR X, \CR Y).\]
\end{itemize}
\end{proposition}


Finally, recall that, similar to the proof of Proposition 2.5 of \cite{AHT}, one can deduce that the kernel and the the cokernel of the unit of adjunction $\delta_X: X \lrt \CE\CR X$ are in $\Add S$. Moreover, $\delta_X$ is a monomorphism if and only if $\Hom_A(S, X)=0$. In the following easy lemma we show that $\delta_X$ is an epimorphism if and only if $\Ext^1_A(S, X)=0$.

\begin{lemma}\label{delta-epi}
With the above notations, $\delta_X$ is an epimorphism if and only if $\Ext^1_A(S,X)=0$.
\end{lemma}

\begin{proof}
Assume that $\delta_X$ is an epimorphism. So we have the short exact sequence
\[0\lrt S_0 \lrt X \st{\delta_X} \lrt \CE\CR X \lrt 0,\]
where $S_0 \in \Add S$. The induced long exact sequence of Ext groups, in view of the facts that $S_0$ is injective and $\CE\CR X \in \Sperp$, implies the result. For the converse, consider the exact sequence
\[0\lrt S_0 \lrt X \st{\delta_X} \lrt \CE\CR X \lrt S_1 \lrt 0,\]
with $S_0, S_1 \in \Add S$. Break out it into two short exact sequences
\[0 \rt S_0 \rt X \rt L \rt 0 \ \ \ {\rm and} \ \ \ 0 \rt L \rt \CE\CR X \rt S_1 \rt 0.\]
The vanishing of $\Ext^1_A(S, X)$, using the first sequence, implies the vanishing of $\Ext^1_A(S, L)$. Now using this, we get from the second sequence that $\Hom_A(S, S_1)=0$. This implies that $S_1=0$. Hence $\delta_X$ is an epimorphism.
\end{proof}

\begin{setup}\label{Setup}
Throughout the paper, $k$ is an algebraically closed field, $B$ is a finite dimensional $k$-algebra and $A$ is the one point extension of $B$ by a fixed finitely generated projective $B$-module $P_0$. We also let $\CR: \Mod A \lrt \Mod B$ denote the restriction functor and $\CE: \Mod B \lrt \Mod A$ denote the extension functor.
\end{setup}

\section{Tilting subcategories}\label{Sec 3-Tilting Sub}
In this section, we recall the definition of tilting subcategories of an abelian category $\SA$, introduced by Beligiannis \cite{B}.

Let us begin by recalling some facts. Let $\SA$ be an abelian category and $\SX$ be a full subcategory of $\SA$. A morphism $\varphi : X \lrt M$, where $M$ is an object of $\SA$, is called a right $\SX$-approximation of $M$ if $X \in \SX$ and every other morphism $\psi: X' \lrt M$, with $X'\in\SX$, factors through $\varphi$. If every object in $\SA$ admits a right $\SX$-approximation, then $\SX$ is called a contravariantly finite subcategory of $\SA$. The notions of a left $\SX$-approximation and a covariantly finite subcategory are defined dually. If $\SX$ is both a contravariantly finite and a covariantly finite subcategory of $\SA$, then it is called a functorially finite subcategory of $\SA$. If for every object $M$ in $\SA$ there exists a monic $0 \lrt M \lrt X$, where $X \in \SX$, then $\SX$ is called a cogenerating subcategory of $\SA$. Generating subcategories are defined dually.

Finally, for an integer $n \geq 0$, we set
\begin{align*}
 \SX^{{\perp}_n} & := \lbrace M \in\SA ~ \vert ~\Ext^n_\SA(\SX, M)=0\rbrace,\\
\SX^{\perp} & :=\lbrace M \in\SA ~ \vert ~\Ext^i_\SA(\SX, M)=0, ~~~\forall~ i\geq 1\rbrace.
\end{align*}
The notions ${}^{{\perp}_n}\SX$ and ${}^{{\perp}}\SX$ are defined dually. Note that $\Ext^0_\SA$ serves for $\Hom_\SA$.

Now we have the necessary background for the following definition of \cite{B}.

\begin{definition}\label{Def(B)}
Let $\SA$ be an abelian category. An additively closed subcategory $\ST$ of $\SA$ is called a tilting subcategory if
\begin{itemize}
\item[$(i)$] $\ST$ is a contravariantly finite subcategory of $\SA$.
\item[$(ii)$] $\ST^\perp=\Fac(\ST).$
\item[$(iii)$] $\ST^{\perp}$ contains a cogenerating subcategory $\SC$ of $\SA$.
\end{itemize}
\end{definition}

Note that if $\SA$ has enough injective objects, then the condition $(iii)$ in the above definition automatically holds true, because in this case, $\ST^{\perp}$ contains the full subcategory of all injective objects of $\SA$ as a cogenerating subcategory.

Recently in \cite[Definition 2.17]{BBOS} the following definition for a tilting subcategory in an abelian category with enough projective objects is studied.

\begin{definition}\label{Def(BBOS)}
Let $\SA$ be an abelian category with enough projectives. An additively closed subcategory $\ST$ of $\SA$ is called a tilting subcategory if
\begin{itemize}
\item[$(i)$] $\ST$ is contravariantly finite in $\SA$.
\item[$(ii)$] $\Ext^1_\SA(T_1, T_2)=0$, for all $T_1, T_2\in\ST$.
\item[$(iii)$]  Every object $T\in\ST$ has projective dimension at most 1.
\item[$(iv)$] For every projective object $P$ in $\SA$, there exists a short exact sequence
\[0\lrt P\lrt T^0\lrt T^1\lrt 0\]
with $T^i\in\ST$.
\end{itemize}
An additively closed full subcategory $\ST$ of $\SA$ which satisfies the conditions $(ii)-(iv)$ of Definition \ref{Def(BBOS)}, is called a weak tilting subcategory.
\end{definition}

Recall that a subcategory $\SX$ of $\SA$ is additively closed if it is closed under taking finite direct sums and direct summands.

\begin{remark}\label{WeakTilting}
An abelian category $\SA$ is called noetherian if every ascending chain of subobjects of an object eventually becomes stationary. By \cite[Proposition 2.20]{BBOS}, if $\SA$ is noetherian, then every weak tilting subcategory of $\SA$ is automatically contravariantly finite. The proof of this is based on the fact that in a noetherian abelian category a notion of trace exists. Although $\Mod R$ is not in general a noetherian category, but trace always exists, see \cite{AF} or page 2 of \cite{CF}. Recall that, for a class $\SX$ of $R$-modules, the trace of $\SX$ in an $R$-module $M$, denoted by $\rm{Tr}_{\SX}(M)$, is defined to be the unique largest submodule of $M$ that belongs to $\Fac(\SX)$. Therefore, it follows, using the same proof as in the proof of Proposition 2.20 of \cite{BBOS}, that any weak tilting subcategory of $\Mod R$ is automatically a tilting subcategory.
\end{remark}

The following theorem is proved by Beligiannis, see Remark 1.9 of \cite{BBOS}. Since there is no published account for this fact, we present a proof here.

\begin{theorem}
Let $\SA=\Mod R$. Let $\ST$ be an additive subcategory of $\SA$. Then $\ST$ is a tilting subcategory in the sense of Definition \ref{Def(B)} if and only if it is a tilting subcategory of $\SA$ in the sense of Definition \ref{Def(BBOS)}.
\end{theorem}

We split out the proof to the following two propositions.

\begin{proposition}\label{(BBOS)Implies(B)}
Let $\ST$ be a tilting subcategory of $\SA=\Mod R$ in the sense of Definition \ref{Def(B)}. Then $\ST$ is a tilting subcategory of $\SA$ in the sense of Definition \ref{Def(BBOS)}.
\end{proposition}

\begin{proof}
We just need to show the validity of the Conditions $(ii)\mbox{-}(iv)$ of \ref{Def(BBOS)}. Obviously, every $T \in \ST$ belongs to $\Fac(\ST)=\ST^\perp$. Hence the Condition $(ii)$ holds true. To show the validity of the Condition $(iii)$, we show that $\Ext^2_R(T, M)=0$, for every $T \in \ST$ and every $M \in \Mod R$. Consider the short exact sequence
\[0\lrt M \lrt E \lrt D \lrt 0\] with $E \in \Inj R$.
This, by applying the functor $\Hom_R(T, - )$, induces the exact sequence
\[0 \lrt \Ext^1_R(T, D) \lrt \Ext^2_R(T, M) \lrt 0\]
of abelian groups. But $E$ is in $\ST^\perp$ and $\ST^\perp=\Fac(\ST)$ is closed under factor modules. So $D \in \ST^\perp$ and hence $\Ext^1_R(T, D)=0$. Therefore $\Ext^2_R(T, M)=0$. To see the validity of the Condition $(iv)$, let $P$ be a projective $R$-module. Consider the short exact sequence
\[0\lrt  P\lrt E\lrt D \lrt 0.\]
Since $E$ is in $\ST^\perp =\Fac(\ST)$, there exists an epimorphism $T\lrt E\lrt 0$, where $T\in\ST$. The map $P \lrt E$ factors through $T$, because $P$ is projective, and hence we get the short exact sequence
\[ 0\lrt P\lrt T \lrt L\lrt 0.\]
Now since $L\in \Fac(\ST)$ and $\ST$ is a contravariantly finite subcategory of $\Mod R$, there exists a short exact sequence
\[ 0\lrt K\lrt T^1\st{h}\lrt L\lrt 0, \]
where $T^1\in \ST$ and $h$ is a right $\ST$-approximation of $L$. By applying the functor $\Hom_\SA(\ST,-)$ to the above short exact sequence and using the fact that $h$ is a right $\ST$-approximation, we get $K\in\ST^{\perp}$.
Consider the pull back diagram
\[\xymatrix{
 & & 0 \ar[d]  & 0 \ar[d] &  \\
 & & K\ar[d] \ar@{=}[r] & K \ar[d] & \\
 0 \ar[r] & P\ar@{=}[d] \ar[r] & T^0 \ar[d] \ar[r] & T^1 \ar[d]\ar[r] & 0 \\
 0 \ar[r] & P\ar[r] & {T} \ar[d]  \ar[r] &L \ar[r]\ar[d] & 0 \\
 & & 0 &  0 & }\]
and use the fact that $\ST^{\perp}$ is closed under extensions, to get that $T^0\in\ST^{\perp}=\Fac(\ST)$.  By the similar argument to that we just applied to $L$, applying this time to $T^0$, we get the exact sequence
\[ 0\lrt Y \lrt T'\lrt T^0\lrt 0,\]
where $T'\in\ST$ and $Y\in\ST^{\perp}$. But $\Ext^1_R(T^0, Y)=0$, so the sequence splits. Therefore $T^0$ is  a summand of $T'$. Hence $T^0\in \ST$ and the short exact sequence
\[ 0\lrt P\lrt T^0\lrt T^1\lrt 0\]
is the desired one.
\end{proof}

\begin{proposition}\label{(B)Implies(BBOS)}
Let $\ST$ be a tilting subcategory of  $\SA=\Mod R$ in the sense of Definition \ref{Def(BBOS)}. Then $\ST$ is a tilting subcategory of $\SA$ in the sense of Definition \ref{Def(B)}.
\end{proposition}

\begin{proof}
We just need to show the validity of the Condition $(ii)$  of Definition \ref{Def(B)}. Let $M\in \Fac(\ST)$. Consider the short exact sequence
\[ 0\lrt K\lrt T\lrt M\lrt 0\]
with $T\in \ST$. Let $T'\in\ST$. By applying the functor $\Hom_{R}(T', -)$ on the above short exact sequence, we get the long exact sequence
\[\Ext^1_R(T', T)\lrt \Ext^1_R(T', M)\lrt \Ext^2_R(T', K)\]
of abelian groups.
By the Conditions $(ii)$ and $(iii)$ of \ref{Def(BBOS)}, we have $\Ext^1_R(T', T)=0$ and  $\Ext^2_R(T', K)=0$. Hence  for all $i\geq 1$, $\Ext^i_R(T', M)=0$. So $M\in \ST^{\perp}$. Conversely, let $M\in{\ST}^{\perp}$ and ${\rm{Tr}}_\ST(M)$  denote
the trace of $\ST$ in $M$, the unique largest submodule of $M$ that belongs to $\Fac(\ST)$. The short exact sequence in the Condition $(iv)$ of \ref{Def(BBOS)}, induces the commutative diagram
\[\xymatrix{ \Hom_R(T^0, {\rm{Tr}}_{T^0}(M)) \ar[r] \ar[d]^{\cong} & \Hom_R(P,  {\rm{Tr}}_{T^0}(M)) \ar[d]^{\beta}\\
\Hom_R(T^0, M) \ar[r]^{\alpha} & \Hom_R(P, M)}\]
where $\beta$ is a monomorphism and, since $T^1\in \ST$, $\alpha$ is an epimorphism. Therefore $\beta$ is an isomorphism and $\rm{Tr}_{T^0}(M)= M$. Hence $M\in\Fac(\ST)$.
\end{proof}

\begin{remark}
Note that Proposition \ref{(BBOS)Implies(B)} is valid with the similar proof for an arbitrary abelian category with enough projective objects. Moreover, Proposition \ref{(B)Implies(BBOS)} holds true for a noetherian abelian category with enough projective objects.
\end{remark}

\section{Extending tilting subcategories}\label{Sec 4-TilSub One point}
In this section we study the behaviour of the tilting subcategories of $\Mod A$ and $\Mod B$ under the restriction functor $\CR$ and the extension functor $\CE$. Such study is motivated by \cite[4.1]{AHT} and provides a wide generalization of their results.

For a subcategory $\ST$ of $\Mod A$, we set $\CR\ST = \{\CR T \ : \  T \in \ST \}$. Similarly, for a subcategory $\ST'$ of $\Mod B$, we set $\CE \ST'= \{\CE M \ : \ M \in \ST'\}.$ Furthermore, we set $\CE\ST' \oplus S := {\Add}\lbrace\CE {\ST'}, {S}\rbrace$ and $\SE\ST' \oplus S := {\add}\lbrace\SE{\ST'}, {S}\rbrace$, where we have fixed $S$ for the simple injective $A$-module $i_*(k)$.

\begin{theorem}\label{Tilting}
\begin{itemize}
\item[$(a)$] Let $\ST$ be a tilting subcategory of $\Mod A$. Then $\ST'=\CR \ST$ is a tilting subcategory of $\Mod B$.
\item[$(b)$] Let $\ST'$ be a tilting subcategory of $\Mod B$. Then $\CE\ST'\oplus S$ is a tilting subcategory of $\Mod A$.
\end{itemize}
\end{theorem}

\begin{proof}
$(a)$ We show that $\ST'$ satisfies the conditions of Definition \ref{Def(B)}. Let $M\in\Fac(\ST')$. There exists a short exact sequence
\[ 0\lrt K\lrt \CR T^1 \lrt M\lrt 0\]
where $T^1\in \ST$. Assume that $\CR T^0 \in \ST'=\CR\ST$. By applying the functor $\Hom_B(\CR T^0,-)$ to the above short exact sequence we get the following exact sequence
\[\Ext^1_B(\CR T^0 , \CR T^1)\lrt \Ext^1_B(\CR T^0 , M)\lrt \Ext^2_B(\CR T^0, K).\]
Since $T^0, T^1\in\ST$ and $\ST$ is a tilting subcategory of $\Mod A$, we have $\Ext^1_A(T^0, T^1)=0$ and $\pd_A T^0\leq 1$. So by the statements $(2)$ and $(4)$ of Proposition \ref{AHT}, we have $\Ext^1_B(\CR T^0 , \CR T^1)=0$ and $\pd_B \CR T^0 \leq 1$. Hence we get $\Ext^1_B(\ST', M)=0$ and therefore $M\in\ST'^{\perp}$.
Now we show the reverse inclusion. Let $M\in\ST'^{\perp}$. By the statement $(2)$ of Proposition \ref{AHT}, $\CE M\in \ST^{\perp}$.
Since $\ST$ is a tilting subcategory of $\Mod A$, $\ST^{\perp}=\Fac(\ST)$. Hence $\CE M\in\Fac(\ST)$ and so there is an epimorphism $f:T\lrt \CE M$ such that $T\in\ST$. By using the restriction and extension sequences we have the following commutative diagram
\[ \xymatrix{ 0\ar[r]  & \CR T \ar[r]\ar[d]^{f'} & T \ar[r] \ar[d]^{f} & S'\ar[r] \ar[d]^{f''}&0\\
 0 \ar[r]& M\ar[r] & \CE M\ar[r] & S''\ar[r]& 0}\]
where $f'=\CR f: \CR T\lrt \CR\CE M\cong M$. Now by the Snake lemma, we have an epimorphism $g: \Ker f''\lrt \Coker f'$. But $\Coker f'=0$, since $\Ker f'' \in \Add S$ and $\Coker f'\in\Mod B$. Therefore, $f'$ is an epimorphism and $M\in\Fac(\CR\ST)$. The last condition is valid because $\Mod B$ has enough injective objects.
Hence $\ST'$ is a weak tilting subcategory of $\Mod B$. Now the result follows from Remark \ref{WeakTilting}.

$(b)$ We proceed as in the proof of part $(a)$. The inclusion $\Fac(\CE \ST'\oplus S)\subseteq(\CE \ST'\oplus S)^{\perp}$ follows from Proposition \ref{AHT}. Now let $X\in (\CE \ST'\oplus S)^{\perp}$. Let $X=X'\oplus S$, where $X'$  does not  have $S$ as a summand. So $\Hom_A(S, X')=0$. On the other hand, since $\Ext^1_A(\CE \ST'\oplus S, X)=0$, we have $\Ext_A^1(S, X')=0$. Therefore, $X'\in \Sperp$. By the statement $(1)$ of Proposition \ref{AHT}, we get $\Ext^1_B(\ST', \CR X')=0$ and so $\CR X'\in\ST'^{\perp}$. But $\ST'$ is a tilting subcategory of $\Mod B$, so $\ST'^{\perp}=\Fac(\ST')$. Hence $\CR X'\in \Fac(\ST')$ and there exists an  epimorphism
\[T'\lrt \CR X'\lrt 0\]
where $T'\in\ST$.
Now since $\CE$ is an exact functor and $X'\in \Sperp$, we get the exact sequence
\[ \CE T'\lrt  \CE \CR X'\cong X'\lrt 0\]
in which we conclude $X'\in \Fac(\CE \ST')$ and $X\in \Fac(\CE \ST'\oplus S)$. The proof now is completed in view of Remark \ref{WeakTilting}.
\end{proof}

Using the same argument as in the proof of the above theorem, similar result could be proved in small mod. So we just state the theorem and skip the proof.

\begin{theorem}\label{Small-AHT}
\begin{itemize}
\item[$(a)$] Let $\ST$ be a tilting subcategory of $\mmod A$. Then $\SR \ST$ is a tilting subcategory of $\mmod B$.
\item[$(b)$] Let $\ST'$ be a tilting subcategory of $\mmod B$. Then $\SE\ST'\oplus S$ is a tilting subcategory of $\mmod A$.
\end{itemize}
\end{theorem}

As an immediate consequence we have the following corollary, that provides another proof for Proposition 4.1 of \cite{AHT}. Recall that a $\La$-module $T$, where $T$ is an artin algebra, is called a tilting module if $\Ext^1_{\La}(T,T)=0$, the projective dimension of $T$, $\pd_{\La}T$ is at most $1$, and there exists a short exact sequence $0 \rt \La \rt T^0 \rt T^1 \rt 0$, such that $T^0$ and $T^1$ are in $\add T$. It is easy to see that if $T$ is a tilting $B$-module , then $\add T$ is a tilting subcategory of $\mmod B$, see \cite{AST}.

\begin{corollary}
\begin{itemize}
\item[$(a)$] Let $T$ be a tilting $A$-module. Then $\SR T$ is a tilting $B$-module.
\item[$(b)$] Let $T'$ be a tilting $B$-module. Then $\SE T' \oplus S$ is a tilting $A$-module.
\end{itemize}
\end{corollary}

\begin{proof}
$(a)$. Let $T$ be a tilting module in $\mmod A$. Then $\add T$ is a tilting subcategory of $\mmod A$. and so by Theorem \ref{Small-AHT}, $\SR \add T$ is a tilting subcategory of $\mmod B$. But clearly $\SR \add T = \add \SR T$. Hence $\SR T$ is a tilting module in $\mmod B$.

The statement $(b)$ proves similarly.
\end{proof}

\begin{remark}\label{addgentilt}
Our proof of the Theorem \ref{Small-AHT} is independent of the main result of \cite{AHT}. However, in the small mod, one can use Proposition 4.1 of \cite{AHT}. Let $\ST$ be a tilting subcategory of $\mmod A$. Since $\prj A = \add A$, by Proposition 3.42 of \cite{R}, there exists a tilting object $T$ in $\ST$ such that $\ST=\add T$ to provide another proof for Theorem \ref{Small-AHT}. In fact, Proposition 4.1 of \cite{AHT} implies that $\SR T$ is a tilting module in $\mmod B$ and hence $\add\SR T$ is a tilting subcategory of $\mmod B$. But $\add\SR T = \SR\add T=\SR\ST$ implies that $\SR\ST$ is a tilting subcategory of $\mmod B$. The same is true for the Statement $(b)$ of the corollary.
\end{remark}

\section{Extending special modules}\label{Sec 5-Special Mod}
In this section we investigate and discuss the restriction and extension of finendo quasi-tilting modules, silting modules and cosilting modules.

Recall that an $R$-module $T$ is Ext-projective, resp. Ext-injective, with respect to a subcategory $\SX$ of $\Mod R$, if $T \in \SX$ and $\Ext^1_R(T, \SX)=0$, resp. $\Ext^1_R(\SX, T)=0$.

\subsection{Silting modules}
Let $\sigma$ be a morphism in $\Prj R$. Let $\SD_{\sigma}$ denote the class of all modules $M$ in $\Mod R$ such that the induced morphism $\Hom_R(\sigma, M)$ is surjective.

\begin{definition}(see \cite{AMV})
Let $T$ be an $R$-module. $T$ is called
\begin{itemize}
  \item finendo, if it is finitely generated over its endomorphism ring.
  \item quasi-tilting, if it is Ext-projective in $\Gen(T)$ and $\Pres(T)=\Gen(T)$.
  \item silting, if there exists a projective presentation $\sigma$ of $T$ such that $\SD_{\sigma}=\Gen(T)$.
\end{itemize}
\end{definition}

If $T$ is both finendo and quasi-tilting it is called finendo quasi-tilting. By \cite[Proposition 3.10]{AMV}, every silting module is finendo quasi-tilting.

For the restriction of a silting module, we need an assumption on the vanishing of Ext.

\begin{theorem}
Let $T$ be a silting $A$-module with respect to a projective presentation $\sigma$. Then the following hold.
\begin{itemize}
\item[$1.$] $\SD_{\CR \sigma}\subseteq \Gen(\CR T)$.
\item[$2.$] If $\Ext^1_A(S, T)=0$, then $\Gen(\CR T) \subseteq \SD_{\CR\sigma}$. In particular, $\CR T$ is a silting $B$-module.
\end{itemize}
\end{theorem}

\begin{proof}
$1.$ Since $T$ is a silting $A$-module, there exists a projective presentation $Q_1\st{\sigma}\lrt Q_0$ of $T$ such that $\SD_{\sigma}=\Gen(T)$. Since the functor $\CR$ is exact and preserves projectives, the induces sequence
\[\CR Q_1\st{\CR \sigma}\lrt \CR Q_0\lrt \CR T\lrt 0\]
is a projective presentation of $\CR T$.

 Let $M\in\SD_{\CR\sigma}$. So there exists an epimorphism
\[\Hom_B(\CR Q_0, M)\lrt \Hom_B(\CR Q_1, M)\lrt 0.\]
Since $(\CR, \CE)$ is an adjiont pair, we get an epimorphism
\[\Hom_A(Q_0, \CE M)\lrt \Hom_A(Q_1, \CE M)\lrt 0\]
which shows that $\CE M\in \SD_{\sigma}=\Gen(T)$. Therefore $\CR\CE M\cong M\in \CR(\Gen(T))$. But since the functor $\CR$ commutes with arbitrary direct sums,  we have $M\in\Gen(\CR T)$. Hence we have proved the statement $1$.

$2$. Since $\Ext^1_A(S, T)=0$, by Lemma \ref{delta-epi}, $\delta_T: T \lrt \CE\CR T$ is an epimorphism. Let $N\in\Gen(\CR T)$. There exists an epimorphism $\bigoplus \CR T\lrt  N\lrt 0$. Thus $\CE N\in \Gen(\CE\CR T)$. Now since $\delta_T$ is an epimorphism, $\CE N\in \Gen(T)$.
But $\Gen(T)=\SD_\sigma$ implies that $\CE N\in\SD_\sigma$. Therefore there exists an epimorphism
\[\Hom_A(Q_0, \CE N)\lrt \Hom_A(Q_1, \CE N)\lrt 0.\]
By using the adjoint pair $(\CR, \CE)$ we have an epimorphism
\[\Hom_B(\CR Q_0, N)\lrt \Hom_B(\CR Q_1, N)\lrt 0\]
which shows that $N\in\SD_{\CR\sigma}$. So $\Gen(\CR T)\subseteq\SD_{\CR\sigma}$ and hence $\SD_{\CR\sigma}=\Gen(\CR T)$. Therefore $\CR T$ is a silting $B$-module.
\end{proof}

\begin{theorem}\label{ExtFinQuasitilting}
Let $B$ and $A$ be as in the Setup \ref{Setup}.
\begin{itemize}
  \item[$1.$] Let $T$ be a finendo quasi-tilting module in $\Mod A$. Then $\CR T$ is a finendo quasi-tilting module in $\Mod B$.
  \item[$2.$]Let $T'$ be a finendo quasi-tilting $B$-module. Then $\CE T' \oplus S$ is a finendo quasi-tilting $A$-module.
\end{itemize}

\end{theorem}

\begin{proof}

$1.$ For the proof of this statement, we use a bijection between finendo quasi-tilting modules and support $\tau$-tilting subcategories which comes in section \ref{Sec 6-tau-tilSub}, after Theorem \ref{Th8.1.10}.

$2.$ First we show that $\CE T'\oplus S$ is $\Ext$-projective in $\Gen(\CE T'\oplus S)$. Let $X\in \Gen(\CE T'\oplus S)$. We may assume that $X=X'\oplus S^{(I)}$, where $S^{(I)}$ is a direct coproduct of $I$ copies of $S$ and $X'$ does not have any copy of $S$ as a summand. Consider an epimorphism
\begin{equation}\label{epi}
\bigoplus (\CE T'\oplus S)\lrt X'\lrt 0.
\end{equation}
By applying the functor $\Hom_A(S, -)$ we get $\Ext^1_A(S, X')=0$. Hence $X'\in \Sperp$ and so $X' \cong \CE\CR X'$.

Now by applying the functor $\CR$ on the \ref{epi}, we get an epimorphism
\begin{equation}
\bigoplus T'\lrt \CR X'\lrt 0,
\end{equation}
which, in turn, implies that $\CR X'\in \Gen(T')$. Therefore $X'=\CE\CR X'\in\CE(\Gen(T'))$. So $X'=\CE M$, for some $M\in\Gen(T')$.
On the other hand, by adjoint property of adjoint pair $(\CR, \CE)$, we have
\[\Ext^1_A(\CE T', \CE M)=\Ext^1_B(T', M).\]
Since $T'$ is a finendo quasi-tilting $B$-module, it is $\Ext$-projective in $\Gen(T)$. So $\Ext^1_B(T', M)=0$. Hence $\Ext^1_A(\CE T', \CE M)=0$. On the other hand, since $A$ is a finite dimensional algebra, $S^{(I)}$ is an injective $A$-module. Thus we conclude that
\[\Ext^1_A(\CE T'\oplus S, X'\oplus S^{(I)})=0,\]
which shows that $\CE T'\oplus S$ is $\Ext$-projective in $\Gen(\CE T'\oplus S)$.

Next, we show that $\Gen(\CE T'\oplus S)\subseteq \Pres(\CE  T'\oplus S)$. Let $X\in \Gen(\CE T'\oplus S)$.  By the similar argument as in the first part of the proof,  we have  $X=X'\oplus S^{(I)}$ such that  $X'\in\CE(\Gen(T'))$. Therefore $\CR X'\in \Gen(T')$. But  since $T'$ is   finendo quasi-tilting,  $\CR X'\in\Pres(T')$.  Thus there exists an exact sequence
\[T'_1\lrt T'_0\lrt \CR X'\lrt 0\]
such that $T'_1, T'_0\in\Add T'$. By applying the exact functor $\CE$ on the above exact sequence and using the fact that $X'\in\Sperp$, we get $X'\in\Pres(\CE T')$. Hence $X=X'\oplus S\in\Pres(\CE T'\oplus S)$.

Finally since $T'$ is finitely generated over its endomorphism ring, we conclude that $\CE T \oplus S$ is finitely generated over its endomorphism ring. The proof is hence complete.
\end{proof}

Let $\La$ be an artin algebra. A module $T\in \mmod \La$ is called a support $\tau$-tilting module \cite{AIR} if $\Ext^1_\SA(T, \Fac(T))=0$ and there exists an exact sequence
\[B \st{f} \lrt T^0 \lrt T^1\lrt 0\]
such that $T^0, T^1\in \add T$ and $f$ is a left $\add T$-approximation, see \cite{J}.

Using previous result, we are able to provide another proof for the main result, i.e. Proposition 3.2, of \cite{Su}.

\begin{corollary}
Let $A=B[P_0]$ as in the Setup \ref{Setup}.
\begin{itemize}
  \item [$(i)$] Let $T$ be a support $\tau$-tilting module of $\mmod A$. Then $\CR T$ is a support $\tau$-tilting module in $\mmod B$.
  \item [$(ii)$] Let $T'$ be a support $\tau$-tilting module of $\mmod B$. Then $\CE T' \oplus S$ is a support $\tau$-tilting module in $\mmod A$.
\end{itemize}
\end{corollary}

\begin{proof}
By the main theorem of \cite{W}, over a finite dimensional algebra, finendo quasi-tilting modules and support $\tau$-tilting modules coincide. Hence the result follows from the Theorem \ref{ExtFinQuasitilting}.
\end{proof}

\begin{remark}
Let $T'$ be a silting $B$-module. So by \cite[Proposition 3.10]{AMV} it is a finendo quasi-tilting module. Hence by the above theorem $\CE T \oplus S$ is a finendo quasi-tilting $A$-module. That is the extension of a silting $B$-module is a finendo quasi-tilting $A$-module. However, it is not clear if the extension of a silting $B$-module is a silting $A$-module.
\end{remark}

\subsection{Cosilting modules}
In this subsection, we study the behaviour of cosilting modules over the one-point extension algebras. Let $\zeta$ be a morphism in $\Inj R$. Let $\SB_{\zeta}$ denote the class of all modules $M$ in $\Mod R$ such that the induced morphism $\Hom_R(M, \zeta)$ is surjective.

\begin{definition}(See \cite[Proposition 1.6]{ATT}, \cite[Definition 2.1]{ZW} and \cite[Definition 3.1]{BP})
Let $T$ be an $R$-module. $T$ is called
\begin{itemize}
  \item cofinendo if there exists a right $\Prod(T)$-approximation of an injective cogenerator $\Mod R$.
  \item quasi-cotilting if it is $\Ext$-injective in $\Cogen(T)$ and $\Copres(T)=\Cogen(T)$.
  \item cosilting if there is an injective copresentation $\zeta$ of $T$ such that $\Cogen(T)=\SB_{\zeta}.$
\end{itemize}
\end{definition}

By \cite[Proposition 2.11]{ZW}, all quasi-cotilting $R$-modules are cofinendo. Moreover, by Theorem 4.18 of \cite{ZW1}, quasi-cotilting modules and cosilting modules are coincide.

\begin{theorem}\label{ExtCosilting}
Let $B$ and $A$ be as in the Setup \ref{Setup}.
\begin{itemize}
  \item[$1.$] Let $C\in\Sperp$ be a cosilting $A$-module. Then $\CR C$ is a cosilting  $B$-module.
  \item[$2.$] Let $C'$ be a cosilting $B$-module. Then $\CE C'\oplus S$ is a cosilting $A$-module.
\end{itemize}
\end{theorem}

\begin{proof}
$1.$ Since cosilting modules and quasi-cotilting modules are coincide, we  show that $\CR C$ is quasi-cotilting. So we need to show that $\CR C$ is $\Ext$-injective in $\Cogen(\CR C)$ and $\Cogen(\CR C)\subseteq \Copres(\CR C)$.
Let $M\in \Cogen(\CR C)$. So $\CE M\in \Cogen(C)$ and therefore  $M\in\CR(\Cogen(C))$.  Hence $M\cong \CR X$, where $X\in \Cogen(C)$. Since $C\in\Sperp$, we have $\Ext^1_B(\CR X, \CR C)\cong \Ext^1_A(\CE\CR X, C)$.
 Since  $C$ is a quasi-cotilting $A$-module and $\CE\CR X =\CE M \in \Cogen(C)$, so $\Ext^1_A(\CE\CR X, C)=0$ and hence $\Ext^1_B(M, \CR C)=0$. Thus  $\CR C$ is $\Ext$-injective in $\Cogen(\CR C)$. Now we show
$\Cogen(\CR C)\subseteq \Copres(\CR C)$. Let $M\in\Cogen(\CR C)$. By the similar argument as in the first part, we have $M\cong \CR X$, where $X\in\Cogen(C)$. So by the assumption $X\in\Copres(C)$. Hence there exists an exact sequence
\[0\lrt X\lrt C^0\lrt C^1\]
such that $C^i\in\Add(C)$. By applying the functor $\CR$ to the above exact sequence we get $\CR X\cong M\in\Copres(\CR C)$, as we wanted.

$2.$ Since $C'$ is a cosilting $B$-module, there exists an injective copresentation\[0\lrt C'\lrt  E^0\st{\zeta}\lrt  E^1\] of $C'$ such that $\Cogen(C')=\SB_{\zeta}.$ By applying the exact functor $\CE$ on the above exact sequence and using the fact that  $\CE$ preserves injective modules,
we get
\[0\lrt \CE C'\oplus S\lrt \CE E^0\oplus S\st{(\CE \zeta\ 0)}\lrt \CE E^1,\]
which is an injective copresentation for $\CE C'\oplus S$. We claim that $\SB_{(\CE\zeta\ o)}=\Cogen(\CE C'\oplus S)$. To prove the claim, first let $X\in\SB_{(\CE\zeta\ 0)}$. We can assume that $S$ is not a direct summand of $X$, since it is clear that $S$ is  in both $\SB_{(\CE\zeta\ o)}$ and $\Cogen(\CE C'\oplus S)$.  Then  the  induced  morphisms
\[\Hom_A(X, \CE E^0\oplus S)\lrt \Hom_A(X, \CE E^1)\]
and hence the morphism
\[\Hom_A(X, \CE E^0)\lrt \Hom_A(X, \CE E^1)\]
are epimorphisms. Therefore since $(\CR, \CE)$ is an adjoint pair, the following sequence
\[ \Hom_B(\CR X, E^0)\lrt \Hom_B(\CR X, E^1)\lrt 0,\]
is exact, which shows that $\CR X\in \SB_{\zeta}=\Cogen(C')$. So there is a monomorphism $0 \lrt \CR X \lrt \prod C'$ that induces the monomorphism $0\lrt \CE\CR X\lrt \CE \prod C'=\prod \CE C'$. Since $S$ is not a summand of $X$, the morphism $\delta_X: X\lrt \CE\CR X$ is a monomorphism. Thus we get the monomorphism $0\lrt X\lrt \prod\CE C'$ which shows that $X\in\Cogen(\CE C'\oplus S)$. Hence $\SB_{(\CE\zeta\ 0)}\subseteq \Cogen(\CE C'\oplus S)$.
For the converse inclusion, let $Y\in \Cogen(\CE C'\oplus S)$. We assume that $Y=Y'\oplus S^{(I)}$ such that $S$ is not a  direct summand of $Y'$.
 Then there exists a monomorphism
\begin{equation*}
0\lrt Y'\lrt \prod(\CE C'\oplus S).
\end{equation*}
By applying the exact functor $\CR$ to the above exact sequence, we get   a monomorphism
\begin{equation*}
0\lrt \CR Y'\lrt \prod C'
\end{equation*}
which shows that $\CR Y'\in \Cogen(C')$. Thus $\CE\CR Y'\in \CE(\Cogen(C'))$. But since $S$ is not a  direct summand of $Y'$, we have monomorphism $0\lrt Y'\lrt \CE\CR Y'$ which shows that $Y'\in\CE(\Cogen(T'))$. So $Y'=\CE M$ where $M\in\Cogen(C')$. Therefore $\CR Y'\in\Cogen(C')=\SB_\zeta$.
By the definition of $\SB_\zeta$, there exists an epimorphism
\[\Hom_B(\CR Y', E^0)\lrt \Hom_B(\CR Y', E^1)\lrt 0.\]
Since $(\CR, \CE)$ is an adjoint pair, we get the epimorphism
\[\Hom_A(Y', \CE E^0)\lrt \Hom_A(Y', \CE E^1)\lrt 0,\]
that, in turn, induces the following epimorphism
\[\Hom_A(Y, \CE E^0)\oplus \Hom_A(Y, S)\lrt \Hom_A(Y, \CE E^1)\lrt 0.\]
Hence $Y\in \SB_{(\CE\zeta\ 0)}$, which shows  $\Cogen(\CE C'\oplus S)\subseteq B_{(\CE\zeta \ 0)}$. This completes the proof of the claim. So $\CE C'\oplus S$ is a cosilting $A$-module.
\end{proof}

\section{$\tau$-tilting subcategories}\label{Sec 6-tau-tilSub}
As an application of the results of the previous section, in view of \cite{AST}, we are able to study the restriction and extension of $\tau$-tilting subcategories. Let us first recall the definition of a $\tau$-tilting subcategory of an abelian category with enough projective objects.

\begin{definition} (see \cite[Definition 1.5]{IJY} and \cite[Definition 2.1]{LZh})\label{DefTauTilting}
Let $\SA$ be an abelian category with enough projective objects. An additively closed subcategory $\ST$ of $\SA$ is called
\begin{itemize}
  \item $\tau$-rigid if for every object $T\in \ST$ there exists an exact sequence \[P_1\st{g}\lrt P_0\lrt T\lrt 0\] such that $P_1$ and $P_0$ are projectives and for every $T'\in \ST$, $\Hom_\SA(g, T')$ is an epimorphism.
  \item weak support $\tau$-tilting if it is a $\tau$-rigid subcategory of $\SA$ and for every projective object $P$ in $\SA$, there exists an exact sequence
  \[P \st{f}{\lrt} T^0 \lrt T^1 \lrt 0\] such that $T^0$ and $T^1$ are in $\ST$ and $f$ is a left $\ST$-approximation of $P$.
  \item support $\tau$-tilting if it is a weak support $\tau$-tilting subcategory which is contravariantly finite in $\SA$.
  \item $\tau$-tilting if it is a support $\tau$-tilting subcategory such that for every projective object $P$ in $\SA$, there exists an exact sequence \[P \st{f}{\lrt} T^0 \lrt T^1 \lrt 0\] such that $T^0$ and $T^1$ are in $\ST$ and $f$ is a non-zero left $\ST$-approximation of $P$.
\end{itemize}
\end{definition}

\begin{theorem}\label{ResTauTilting}
Let $\ST$ be a support $\tau$-tilting subcategory of $\Mod A$. Then $\CR\ST$ is a support $\tau$-tilting subcategory of $\Mod B$.
\end{theorem}

\begin{proof}
First we show that $\CR\ST$ is a $\tau$-rigid subcategory of $\Mod B$. Since $\ST$ is a $\tau$-rigid subcategory of $\Mod A$, for every object $T\in\ST$, there exists a projective presentation
\[Q_1 \st{g} \lrt Q_0 \lrt T \lrt 0\]
of $T$ such that the induced morphism $\Hom_A(g, \ST)$ is an epimorphism. We can deduce that $Q_1$ does not contain $\tilde{P}$ as a summand. Otherwise, $S$ should appear as a summand of $Q_0$, i.e. should be projective, which contradicts with the fact that $P_0 \neq 0$. Since $\CR$ is an exact functor and preserves projective modules, the exact sequence
\[\CR Q_1\st{\CR g}\lrt \CR Q_0\lrt \CR T\lrt0\] is a projective presentation of $\CR T$.
Now in order to show the result, we have to show that $\Hom_B(\CR g, \CR\ST)$ is an epimorphism.
Let $\CR T'\in\CR\ST$  and $f\in\Hom_B(\CR Q_1, \CR T')$. Since $\CR Q_1\cong Q_1$, the morphism $f$ induces morphism $\tilde f\in \Hom_B(Q_1, T')$ such that $\tilde{f}=if$ where $i: \CR T'\lrt T'$ is the inclusion. But, $\tau$-rigidity of $\ST$ implies that there exists a morphism $h: Q_0\lrt T'$ such that $\tilde{f}=hg$. Therefore $\CR\tilde{f}=\CR h\CR g$. Hence $f=\CR h\CR g$.

Since $\ST$ is a $\tau$-tilting subcategory of $\Mod A$, for every projective $B$-module $P$, there exists an exact sequences
\[P\st{f}\lrt T^0\lrt T^1\lrt 0 \]
such that $T^0, T^1\in\ST$ and $f$ is a   left $\ST$-approximation of $P$.  By applying the exact functor $\CR$ on the above exact sequence we have the following exact sequence
\[\CR P\st{\CR f}\lrt \CR T^0\lrt \CR T^1\lrt 0.\]
Here  $\CR P\cong P\in\Prj B$, so in order to show that this is a desired exact sequence, we need to show that $\CR f$ is a   left $\CR\ST$-approximation. Let $h: \CR P\lrt \CR U$ be a morphism such that $U\in\ST$. Then the  morphism $h$ induces morphism $\tilde h\in \Hom_B(P, U)$ such that $\tilde{h}=ih$ where $i: \CR U\lrt U$ is the natural inclusion.  Since $f$ is a left $\ST$-approximation, there exists  morphism $g: T^0\lrt U$ such that $gf=\tilde{h}$. Therefore $\CR g\CR f=\CR \tilde{h}=h$.

To complete the proof, we just need to  show that, $\ST'$ is a contravariantly finite subcategory of $\Mod B$. Let $M\in \Mod B$.  Since $\ST$ is a contravariantly finite  subcategory of $\Mod A$, there is a right $\ST$-approximation $T\lrt \CE M$, with $T\in\ST$, for $\CE M\in \Mod A$. By applying the  functor $\CR$ we get a right $\ST'$-approximation $\CR T\lrt M$ for $M$. Hence $\ST'$ is a contravariantly finite subcategory of $\Mod B$.
\end{proof}

We say that two support $\tau$-tilting subcategories $\ST$ and $\ST'$ are equivalent if $\Fac(\ST)=\Fac(\ST')$. Moreover, two quasi-tilting $B$-modules $T$ and $T'$ are called equivalent if $\Add T =\Add T'$. We need the following theorem.

\begin{theorem}(see \cite[Theorem 8.1.10]{AST})\label{Th8.1.10}
Let $R$ be a ring. There is a bijection between the collection of equivalence classes of $\tau$-tilting subcategories of $\Mod R$ and the collection of equivalence classes of finendo quasi-tilting $R$-modules. Based on this bijection, $\Add T$ for any finendo quasi-tilting $R$-module is a support $\tau$-tilting subcategory of $\Mod R$.
\end{theorem}

Now we can prove the  part $1$ of Theorem \ref{ExtFinQuasitilting}.

\vspace{0.2cm}

\noindent \emph{Proof of Theorem \ref{ExtFinQuasitilting}.$1$}.
Since $T$ is a finendo quasi-tilting module, by Theorem \ref{Th8.1.10}, $\Add T$ is a support $\tau$-tilting subcategory of $\Mod A$. Hence by Theorem \ref{ResTauTilting}, $\CR \Add T$ is a $\tau$-tilting subcategory of $\Mod B$. It is obvious that $\CR \Add T=\Add \CR T$. Now the result follows from Proposition 8.1.7 of \cite{AST}.

\begin{theorem}\label{ExtTauTilting}
Let $\ST'$ be a support $\tau$-tilting subcategory of $\Mod B$. Then $\CE \ST'\oplus S$ is a support $\tau$-tilting subcategory of $\Mod A$.
\end{theorem}

\begin{proof}
By Theorem \ref{Th8.1.10}, there is a finendo quasi-tilting $B$-module $T'$ such that $\Fac(\ST')=\Gen(T')$. Theorem \ref{ExtFinQuasitilting}.2, implies that $\CE T'\oplus S$ is a finendo quasi-tilting $A$-module. Theorem 8.1.5 of \cite{AST}, implies that $\Add(\CE T'\oplus S)$ is a support $\tau$-tilting subcategory of $\Mod A$. To conclude the result, it is enough to note that $\Fac(\Add(\CE T'\oplus S))=\Fac(\CE\CT' \oplus S)$. This follows using the fact that $\Fac(\Add(\CE T'\oplus S))= \Gen(\CE T'\oplus S)$.
\end{proof}

\begin{corollary}
\begin{itemize}
  \item [$(i)$] Let $\ST$ be a support $\tau$-tilting subcategory of $\mmod A$. Then $\SR\ST$ is a support $\tau$-tilting subcategory of $\mmod B$.
  \item [$(ii)$] Let $\ST'$ be a support $\tau$-tilting subcategory of $\mmod B$. Then $\SE\ST' \oplus S$ is a support $\tau$-tilting subcategory of $\mmod A$.
\end{itemize}
\end{corollary}

\begin{proof}
The proof follows using similar arguments as in the proofs of Theorems \ref{ResTauTilting} and \ref{ExtTauTilting}.
\end{proof}

\section{$(\tau\mbox{-})$Cotorsion torsion triples}\label{Sec 7-CotTorTriples}
In this section, using bijections between tilting subcategories, resp. support $\tau$-tilting subcategories, and cotorsion torsion triples, resp. $\tau$-cotorsion torsion triples, we study the restriction and extension of these triples. We start by recalling the bijections.

Let $\SA$ be an abelian category with enough projective objects. A pair $(\ST,\SF)$ of full subcategories of $\SA$ is called a torsion pair  if $\Hom_{\SA}(\ST,\SF)=0$ and for every $A \in \SA$ there exists a short exact sequence $0 \lrt tA \lrt A \lrt fA \lrt 0$ such that $tA \in \ST$ and $fA \in \SF$.

A pair $(\SC, \SD)$ of full subcategories of $\SA$ is called a cotorsion pair if $\SC = {}^{\perp_1}\SD$ and $\SD = \SC^{\perp_1}$ and
for every object $A \in \SA$, there exist short exact sequences
\[0 \lrt D \lrt C \lrt A \lrt 0 \ \ {\rm and} \ \ 0 \lrt A \lrt D' \lrt C' \rt 0,\]
where $C$ and $C'$ are in $\SC$ and $D$ and $D'$ are in $\SD$.

A triple $(\SC, \ST, \SF)$  of full subcategories in $\SA$ is called {cotorsion torsion triple}, if $(\SC, \ST)$ is a cotorsion pair and $(\ST, \SF)$ is a torsion pair.

Following theorem provides a tight connection between tilting subcategories of an abelian category with enough projective objects and its cotorsion torsion triples.

\begin{theorem}(see \cite[Theorem 2.29]{BBOS})\label{2.29}
Let $\SA$ be an abelian category with enough projective objects. Then there is a bijection between the collection of all tilting subcategories of $\SA$ and the collection of all cotorsion torsion triples of $\SA$. Based on this bijection, a tilting subcategory $\ST$ maps to the triple
$({}^{\perp_1}{(\rm{Fac}(\ST))}, \rm{Fac}(\ST), \ST^{\perp_0})$ and a cotorsion torsion triples $(\SC, \ST, \SF)$ maps to $\SC \cap\ST$.
\end{theorem}

As a counterpart to a cotorsion torsion triple we have the notion of $\tau$-cotorsion torsion triple. A pair of full subcategories $(\SC, \SD)$ of $\SA$ is called a $\tau$-cotorsion pair \cite{AST} if $\SC\cap\SD$ is a contravariantly finite subcategory of $\SA$, $\SC = {}^{\perp_1}\SD$ and for every projective object $P \in \SA$, there exists an exact sequence
\[  P \st{f}\lrt {D}\lrt C \lrt 0,\]
where $D \in \SC\cap\SD$, $C \in \SC$ and $f$ is a left $\SD$-approximation. A triple $(\SC, \ST, \SF)$ of full subcategories in $\SA$ is called a $\tau$-cotorsion torsion triple, if $(\SC, \ST)$ is a $\tau$-cotorsion pair and $(\ST, \SF)$ is a torsion pair.

There is also a connection between $\tau$-tilting subcategories of an abelian category with enough projective objects and its $\tau$-cotorsion torsion triples.

\begin{theorem}(see \cite[Theorem 5.7]{AST})\label{Bijection}
Let $\SA$ be an abelian category with enough projective objects. Then there is a bijection between the collection of all support $\tau$-tilting subcategories of $\SA$ and the collection of all $\tau$-cotorsion torsion triples of $\SA$. Based on this bijection, a support $\tau$-tilting subcategory $\ST$ maps to the triple
$({}^{\perp_1}{(\rm{Fac}(\ST))}, \rm{Fac}(\ST), \ST^{\perp_0})$ and a $\tau$-cotorsion torsion triples $(\SC, \ST, \SF)$ maps to $\SC\cap\ST$.
\end{theorem}

\subsection{Cotorsion torsion triples}
Now we have the necessary background for our first theorem.

\begin{theorem}\label{ResCotTor}
Let $(\SC, \ST, \SF)$ be a cotorsion torsion triple in $\Mod A$. Then the triple
\[({}^{{\perp_1}}(\Fac(\CR(\SC\cap \ST))), \Fac(\CR(\SC\cap\ST)), {(\CR(\SC\cap\ST))}^{\perp_0})\]
is a cotorsion torsion triple in $\Mod B$. If furthermore, $\ST\subseteq \Sperp$ then $(\CR(\SC), \CR(\ST), \CR(\SF\cap \Sperp))$ is also a cotorsion torsion triple in $\Mod B$ and these two are the same.
\end{theorem}

\begin{proof}
By Theorem \ref{2.29}, $\SC\cap \ST$ is a tilting subcategory of $\Mod A$. Hence by Theorem \ref{Tilting}.$a$, $\CR(\SC\cap\ST)$ is a tilting subcategory of $\Mod B$. Therefore, by using the bijection of Theorem \ref{2.29}, we may conclude the first statement, that is,
\[({}^{{\perp_1}}(\Fac (\CR(\SC\cap \ST))), \Fac (\CR(\SC\cap\ST)), {(\CR(\SC\cap\ST))}^{\perp_0})\]
is a cotorsion torsion triple in $\Mod B$.
For the second statement, first we show  that
\[ \CR(\Fac(\SC\cap\ST))=\Fac (\CR(\SC\cap \ST)).\]
Let $M\in \CR(\Fac(\SC\cap\ST))$. Then $M\cong \CR X$, where $X\in \Fac(\SC\cap\ST)$. Since $X\in \Fac(\SC\cap\ST)$, there exists  an epimorphism
$Y\lrt X\lrt 0$, where $Y\in \SC\cap\ST$. By applying the exact functor $\CR$ on this epimorphism we get $\CR Y\lrt \CR X\cong M\lrt 0$. Therefore $M\in \Fac(\CR(\SC\cap \ST))$. So we have the inclusion
$\CR(\Fac(\SC\cap\ST))\subseteq \Fac (\CR(\SC\cap \ST))$. For the reverse inclusion, let $N\in \Fac (\CR(\SC\cap \ST))$ and consider the epimorphism \[\CR X\lrt N\lrt 0,\]
where $X\in \SC\cap\ST$. Since $\CE$ is an exact functor and $\ST\subseteq \Sperp$, we conclude that $N\in \CR(\Fac(\SC\cap\ST))$. Hence we have the equality $\Fac (\CR(\SC\cap \ST)) = \CR(\Fac(\SC\cap\ST))$. On the other hand,  by Theorem \ref{2.29}, $\SC\cap\ST$ is a tilting subcategory of $\Mod A$ and $\Fac(\SC\cap\ST)=\ST$. So \[\CR(\Fac(\SC\cap\ST))=\CR \ST =\Fac (\CR(\SC\cap \ST)).\]

Now we show that ${}^{{\perp_1}}(\Fac (\CR(\SC\cap \ST)))={}^{{\perp_1}}\CR \ST=\CR \SC$. We have already  seen the first equality. To see the second one, we use the fact that $(\SC, \ST)$ is a cotorsion pair. Since $\SC={}^{\perp_1}{\ST}$, we have $\CR \SC=\CR({}^{\perp_1}{\ST})$. To complete this part of the proof, it is enough to show that  $\CR({}^{\perp_1}{\ST})={}^{{\perp_1}}\CR\ST$. To see this, let $M\in \CR({}^{\perp_1}{\ST})$. Then  $M\cong \CR X $, where $X\in {}^{\perp_1}{\ST}$. Since $X\in {}^{\perp_1}{\ST}$, we have $\Ext^1_A(X, \ST)=0$. So by the statement $(3)$ of Proposition \ref{AHT}, we get $\Ext^1_B(\CR X, \CR\ST)=0$. Therefore $M\cong \CR X\in {}^{{\perp_1}}\CR\ST$.
Now let $N\in {}^{{\perp_1}}(\CR\ST)$. Then $\Ext^1_B(N, \CR\ST)=0$. Since $\ST \subseteq \Sperp $, by the statement $(1)$ of Proposition \ref{AHT}, we obtain $\Ext^1_A(\SE N,\ST)=0$. Hence $\CE N\in {}^{\perp_1}{\ST}$ and $N\cong \CR\CE N\in\CR({}^{\perp_1}{\ST})$.

To complete the proof, we just need to show that ${(\CR(\SC\cap\ST))}^{\perp_0}=\CR(\SF\cap \Sperp)$. To this end, since $\SC\cap\ST$ is a tilting subcategory in $\Mod A$ and $(\SC\cap \ST)^{\perp_0}=\SF$, we show that ${(\CR(\SC\cap\ST))}^{\perp_0}=\CR((\SC\cap\ST)^{\perp_0}\cap \Sperp)$. Let $M\in {(\CR(\SC\cap\ST))}^{\perp_0}$, then $\Hom_B(\CR(\SC\cap \ST), M)=0$. But $(\CR, \CE)$ is an adjoint pair and we get $\Hom_A(\SC\cap\ST, \CE M)=0$.
So $\CE M\in(\SC\cap\ST)^{\perp_0}$. Also $\CE M \in \Sperp$ implies that $\CE M\in(\SC\cap\ST)^{\perp_0} \cap \Sperp$.  Hence $M\cong \CR\CE M \in \CR((\SC\cap\ST)^{\perp_0}\cap \Sperp)$. For the reverse inclusion, let $N\in \CR((\SC\cap\ST)^{\perp_0}\cap \Sperp)$. So $N\cong \CR X$, where $X\in (\SC\cap\ST)^{\perp_0}\cap \Sperp$.
Therefore  we have $\Hom_A(\SC\cap\ST, \CE\CR X)=0$. Hence $\Hom_B(\CR(\SC\cap\ST), \CR X)=0$. So $N\cong \CR X \in {(\CR(\SC\cap\ST))}^{\perp_0}$.

\end{proof}

\begin{theorem}\label{ExtCotTor}
Let $(\SC', \ST', \SF')$ be a cotorsion torsion triple in $\Mod B$. Then the triple
\[ ({}^{{\perp_1}}(\Fac (\CE(\SC'\cap \ST'))), \Fac (\CE(\SC'\cap\ST')\oplus S), {(\CE(\SC'\cap\ST')\oplus S})^{\perp_0})\]
is a cotorsion torsion triple in $\Mod A$. Moreover, the triple $(\CE\SC', \CE\ST', \CE\SF')$ is a cotorsion torsion triple in $\Sperp$ which is equal to the triple
\[({}^{{\perp_1}}(\Fac (\CE(\SC'\cap \ST')))\cap \Sperp, \Fac (\CE(\SC'\cap\ST'))\cap \Sperp, {(\CE(\SC'\cap\ST')})^{\perp_0})\cap \Sperp).\]
\end{theorem}

\begin{proof}

By  Theorems \ref{2.29} and \ref{Tilting}.$b$ we have
\[({}^{{\perp_1}}(\Fac (\CE(\SC'\cap \ST')\oplus S)), \Fac (\CE(\SC'\cap\ST')\oplus S), {(\CE(\SC'\cap\ST')\oplus S})^{\perp_0})\]
is a cotorsion torsion triple in $\Mod A$.  But it is not hard  to see that
\[{}^{{\perp_1}}(\Fac (\CE(\SC'\cap\ST') \oplus S))= {}^{{\perp_1}}(\Fac (\CE(\SC'\cap\ST')) )\oplus {}^{{\perp_1}} \Fac(S),\] for instance see the proof of part $(a)$ of  Proposition 3.2 of \cite{Su}.
On the other hand, we have  ${}^{{\perp_1}}\Fac(S)=0$; indeed, if $X\in \Fac(S)$, then $X\cong S^{k}$, with $k\geq 0$. Since $S$ is an injective module, we have  ${}^{{\perp_1}}\Fac(S)=0$. Hence we have the first statement.
Now we show the second statement. It is plain that the triple $(\CE\SC', \CE\ST', \CE\SF')$ is a cotorsion torsion triple in $\Sperp$. So to complete the proof, it remains to  show that these two triples are the same.
First we show $\Fac (\CE(\SC'\cap\ST'))\cap \Sperp=\CE\ST'$. By Theorem \ref{2.29}, $\SC'\cap\ST'$ is a tilting subcategory of $\Mod B$ and $\Fac(\SC'\cap\ST')=\ST'$. So to show the equality, it is enough to note that $\Fac (\CE(\SC'\cap\ST'))\cap\Sperp=\CE(\Fac(\SC'\cap \ST'))$. Now we consider  \[{}^{{\perp_1}}(\Fac (\CE(\SC'\cap \ST')))\cap \Sperp={}^{{\perp_1}}(\CE(\Fac(\SC'\cap \ST')))\cap\Sperp={}^{{\perp_1}} \CE\ST' \cap \Sperp.\]
Let $X\in {}^{{\perp_1}} \CE \ST'\cap \Sperp$ then $\Ext^1_A(X, \CE\ST')=0$. By the statement $(2)$ of Proposition \ref{AHT}, we have $\Ext^1_B(\CR X, \ST')=0$. Therefore, $\CR X\in {}^{\perp_1}\ST'$. Since $X\in\Sperp$, we have $X\cong\CE\CR X\in \CE({}^{\perp_1}\ST')$. Hence we show the inclusion ${}^{{\perp_1}} \CE \ST' \cap \Sperp\subseteq\CE({}^{\perp_1}\ST')$. Conversely, let $Y\in \CE({}^{\perp_1}\ST')$. Then $Y\cong \CE M $, where $M\in{}^{\perp_1} \ST' $.   Since $M\in{}^{\perp_1} \ST' $, we have $\Ext^1_B(M, \ST')=0$. Hence by the statement $(2)$ of Proposition \ref{AHT}, we get $\Ext^1_A(\CE M, \CE\ST')=0$. So $Y\cong\CE M\in{}^{\perp_1}  \CE\ST' $. We have already shown that  ${}^{{\perp_1}} \CE\ST'\cap \Sperp=\CE({}^{\perp_1}\ST')$.
 Now we note that, since $(\SC', \ST')$ is a cotorsion pair, then $\SC'={}^{\perp_1}\ST'$ and hence $\CE({}^{\perp_1}\ST')=\CE\SC'$. So we get
\[{}^{{\perp_1}}(\Fac (\CE(\SC'\cap\ST')))\cap \Sperp= \CE \SC'.\]
For completing the proof, it remains to show that ${(\CE(\SC'\cap\ST')})^{\perp_0}\cap \Sperp=\CE\SF'$.
First we show that ${(\CE(\SC'\cap\ST'))}^{\perp_0}\cap \Sperp=\CE{(\SC'\cap\ST')}^{\perp_0}$. Let $X\in {(\CE(\SC'\cap\ST'))}^{\perp_0}\cap \Sperp$, then $\Hom_A(\CE(\SC'\cap\ST'), X)=0$. Since $X\in\Sperp$, by the statement $(1)$ of Proposition \ref{AHT}, we have $\Hom_B(\SC'\cap\ST', \CR X)=0$ and hence $\CR X\in {(\SC'\cap\ST')}^{\perp_0}$. Therefore $X\cong \CE\CR X\in \CE({(\SC'\cap\ST')}^{\perp_0})$.
Conversely, let $Y\in \CE{(\SC'\cap\ST')}^{\perp_0}$. Then  $Y\in\Sperp $ and $Y\cong \CE M$, where $M\in {(\SC'\cap\ST')}^{\perp_0}$. Since $M\in {(\SC'\cap\ST')}^{\perp_0}$, we have $\Hom_B(\SC'\cap\ST', M)=0$. By the statement $(1)$ of Proposition \ref{AHT}, $\Hom_A(\CE(\SC'\cap \ST'), \CE M)=0$. Therefore, $Y\cong\CE M\in {(\CE(\SC'\cap\ST'))}^{\perp_0}$. On the other hand, by Theorem \ref{2.29}, ${(\SC'\cap\ST')}^{\perp_0}=\SF'$.
So \[ {(\CE(\SC'\cap\ST'))}^{\perp_0}\cap \Sperp=\CE{(\SC'\cap\ST')}^{\perp_0}=\CE\SF'.\]
Hence the proof is complete.
\end{proof}

\subsection{$\tau$-cotorsion torsion triples}
In this subsection we study the restriction and extension of $\tau$-cotorsion torsion triples, using the bijection mentioned in Theorem \ref{Bijection}.

\begin{theorem}
Let $(\SC, \ST, \SF)$ be a $\tau$-cotorsion torsion triple in $\Mod A$. Then the triple
\[({}^{{\perp_1}}(\Fac(\CR(\SC\cap \ST))), \Fac(\CR(\SC\cap\ST)), {(\CR(\SC\cap\ST))}^{\perp_0})\]
is a $\tau$-cotorsion torsion triple in $\Mod B$. If furthermore, $\ST\subseteq \Sperp$ then $(\CR\SC, \CR\ST, \CR(\SF\cap \Sperp))$ is also a $\tau$-cotorsion torsion triple in $\Mod B$ and these two are the same.
\end{theorem}

\begin{proof}
By Theorems \ref{Bijection}, $\SC\cap\ST$ is a $\tau$-tilting subcategory of $\Mod A$. So by Theorem \ref{ResTauTilting}, $\CR(\SC\cap\ST)$ is a support $\tau$-tilting subcategory of $\Mod B$. So the first statement follows from Theorem \ref{Bijection}. The second part, follows by the similar argument as in the second part of Theorem \ref{ResCotTor}.
\end{proof}

\begin{theorem}
Let $(\SC', \ST', \SF')$ be a $\tau$-triple in $\Mod B$. Then the triple
\[ ({}^{{\perp_1}}(\Fac (\CE(\SC'\cap \ST'))), \Fac (\CE(\SC'\cap\ST')\oplus S), {(\CE(\SC'\cap\ST')\oplus S})^{\perp_0})\]
is a $\tau$-triple in $\Mod A$. Moreover, the triple
\[({}^{{\perp_1}}(\Fac(\CE(\SC' \cap \ST'))) \cap \Sperp, \Fac (\CE(\SC'\cap\ST'))\cap \Sperp, {(\CE(\SC'\cap\ST')})^{\perp_0})\cap \Sperp)\] is a $\tau$-triple in $\Sperp$ and is equal to the triple $(\CE\SC', \SE\ST', \SE\SF')$.
\end{theorem}

\begin{proof}
By Theorems \ref{Bijection}, $\SC' \cap \ST'$ is a $\tau$-tilting subcategory of $\Mod B$. Hence by Theorem \ref{ExtTauTilting}, $\CE(\SC' \cap \ST') \oplus S$ is a $\tau$-tilting subcategory of $\Mod A$. Hence the first statement follows by another use of Theorem \ref{Bijection}. The second part, follows by the similar argument as in the proof of Theorem \ref{ExtCotTor}.
\end{proof}

\section*{Acknowledgments}
This work was partly done during a visit of the first author to the Institut des Hautes \'{E}tudes Scientifiques (IHES), Paris, France. The first and fourth authors would like to thank the support and excellent atmosphere at IHES. The work of the first and the third author is based upon research funded by Iran National Science Foundation (INSF) under project No. 4001480. The fourth author is supported by the European Union’s Horizon 2020 research and innovation programme under the Marie Sklodowska-Curie grant agreement No 893654.


\begin{thebibliography}{29}
\bibitem{AIR}{\sc T. Adachi, O. Iyama and I. Reiten,}{\sl $\tau$-tilting theory,} { Compos. Math.} {\bf 150}(3) (2014), 415-452.

\bibitem{AF}{\sc F. W. Anderson, K. R. Fuller,} {\sl Rings and categories of modules,} Springer, New York, 1992.

\bibitem{AH}{\sc L. Angeleri H\"{u}gel and M. Hrbek,} {\sl Silting modules over commutative rings,} Int. Math. Res. Notices (IMRN) {\bf 13} (2017), 4131-4151.

\bibitem{AMV}{\sc L. Angeleri H\"{u}gel, F. Marks and J. Vit\'{o}ria,} {\sl Silting modules,} Int. Math. Res. Notices (IMRN) {\bf 4} (2016), 1251-1284.

\bibitem{ATT} {\sc L. Angeleri H\"{u}gel, A. Tonolo and J. Trlifaj,}  {\sl Tilting preenvelopes and cotilting precovers,}
Algebr. Represent. Theor.  4(2) (2001), 155-170.

\bibitem{AST} {\sc J. Asadollahi, S. Sadeghi, H. Treffinger}, {\sl On  $\tau$-tilting subcategories}, arXiv:2207.00457v1 [math.RT].

\bibitem{AHT} {\sc I. Assem,  D. Happel and S. Trepode,} {\sl Extending tilting modules to one-point extensions by projectives}, Comm. Algebra {\bf 35}(10) (2007), 2983-3006.

\bibitem{BBOS} {\sc U. Bauer, M. B. Botnan, S. Oppermann and J. Steen,} {\sl  Cotorsion torsion triples and the representation theory of filtered hierarchical clustering,}  Adv. Math. {\bf 369} (2020), 107171.

\bibitem{B} {\sc A. Beligiannis,} {\sl Tilting theory in Abelian categories and related homological and homotopical structures,} 2010, unpublished manuscript.

\bibitem{BP} {\sc S. Breaz and  F. Pop,} {\sl Cosilting modules}, Algebr. Represent. Theor. (2017) 20:1305-1321.

\bibitem{CF} {\sc R. R. Colby and K. R. Fuller,} {\sl Equivalnce and duality for module categories (with tilting and cotilting for rings),} Cambridge University Press, (2004).

\bibitem{IJY} {\sc O. Iyama, P. J\o rgensen and D.Yang}, {\sl Intermediate co-t-structures, two-term silting objects, $\tau$-tilting modules, and torsion classes}, Algebra Number Theory, {\bf 8}(10) (2014), 2413-2431.

\bibitem{J} {\sc G. Jasso,} {\sc Reduction of $\tau$-tilting modules and torsion pairs,} Int. Math. Res. Not. IMRN {\bf 16} (2015), 7190-7237.

\bibitem{KV} {\sc B. Keller, D. Vossieck,} {\sl Aisles in derived categories}, Bull. Soc. Math. Belg. Sér. A {\bf 40}(2) (1988), 239-253.

\bibitem{LZh} {\sc Y. Liu and P. Zhou,} {\sl $\tau$-tilting theory in abelian categories}, Proc. Amer. Math. Soc. {\bf 398} (2014), 63-110.

\bibitem{P} {\sc C. Psaroudakis,} {\sl Homological theory of recollements of abelian categories}, J. Algebra {\bf 398} (2014), 63-110.

\bibitem{R} {\sc E. S. Rundsveen,} {\sl Torsion, cotorsion and tilting in abelian categories,} Master’s thesis in Mathemathical Sciences, NTNU, (2021). Available at: https://ntnuopen.ntnu.no/ntnu-xmlui/bitstream/handle/11250/2778405/no.ntnu:inspera:77742017:16428509.pdf?sequence=1

\bibitem{Su} {\sc P. Suarez,} { \sl $\tau$-tilting modules over one-point extensions by a projective module}, Algebr. Represent. Theor. {\bf 21}(4) (2018), 769-786.

\bibitem{W} {\sc  J. Wei,} {\sl $\tau$-tilting theory and $*$-modules, } J. Algebra {\bf{414 }}(2014), 1-5.

\bibitem{ZW1} {\sc P. Zhang and J. Wei,} {\sl Cosilting complexes and AIR-cotilting modules}, J. Algebra, {\bf 491} (2017), 1-31.

\bibitem{ZW} {\sc P. Zhang and J. Wei,} {\sl Quasi-cotilting modules and torsion-free classes}, J. Algebra Appl. 18, No. 11 (2019) 1950214.
\end{thebibliography}
\end{document}